\documentclass[11pt,a4paper]{article}

\usepackage{makeidx}
\makeindex            

\usepackage{amsmath}  
\usepackage{amssymb}  
\usepackage{amsthm}   
\usepackage{amsfonts}
\newcommand{\N}{\mathbb{N}}

\usepackage{mathrsfs}  

\newcommand{\Ann}[0]{\operatorname{Ann}}
\newcommand{\Hom}[0]{\operatorname{Hom}}
\newcommand{\id}[0]{\operatorname{id}}
\newcommand{\im}[0]{\operatorname{im}}

\newcommand{\Supp}[0]{\operatorname{Supp}}

\newcommand{\Spec}[0]{\operatorname{Spec}}
\newcommand{\End}[0]{\operatorname{End}}
\newcommand{\invlim}[0]{\operatorname{invlim}}
\newcommand{\dirlim}[0]{\operatorname{dirlim}}
\newcommand{\Ext}[0]{\operatorname{Ext}}
\newcommand{\height}[0]{\operatorname{height}}

\newtheorem{satz}{Theorem}[section]
\newtheorem{lemma}[satz]{Lemma}
\newtheorem{coro}[satz]{Corollary}

\newtheorem{definition}[satz]{Definition}
\newtheorem{remark}[satz]{Remark}

\newtheorem{question}[satz]{Question}

\newtheorem*{satzo}{Theorem}

\author{M. Hellus and J. St\"uckrad}
\title{On endomorphism rings of local cohomology modules}
\date{\today}   

\begin{document}

\maketitle

\begin{abstract}
Let $R$ be a local complete ring. For an $R$-module $M$ the
canonical ring map $R\to \End _R(M)$ is in general neither injective
nor surjective; we show that it is bijective for every local
cohomology module $M:=H^h_I(R)$ if $H^l_I(R)=0$ for every $l\neq
h(=\height (I))$ ($I$ an ideal of $R$); furthermore the same holds
for the Matlis dual of such a module. As an application we prove new
criteria for an ideal to be a set-theoretic complete intersection.

\end{abstract}
\setcounter{section}{0}
\section{Introduction}
For an ideal $I$ of a local ring $(R,m)$ we denote the $n$-th local
cohomology functor with support in $I$ by $H^n_I$ and the
(contravariant) Matlis dual functor by $D$, i.~e. $D(M)=\Hom
_R(M,E)$ for any $R$-module $M$, where $E:=E_R(R/M)$ is a fixed
$R$-injective hull of the residue field $R/m$.

Let $R$ be a (always commutative, unitary) ring and $M$ an
$R$-module. There is the canonical map $\mu _M:R\to \End _R(M)$ that
maps $r\in R$ to multiplication by $r$ on $M$; it is a homomorphism
of (associative) $R$-algebras. In general, $\mu _M$ is neither
injective nor surjective. In section 2 we will show that, if $R$ is
noetherian local complete and $I$ an ideal of $R$ such that
$H^l_I(R)=0$ for every $l\neq h(=\height (I))$, then $\mu
_{H^h_I(R)}$ is bijective. In particular, the endomorphism ring of
the $R$-module $H^h_I(R)$ is commutative and $\Ann _R(H^h_I(R))=0$.

The proof of this result uses a generalization of theorem 3.2 from
\cite{hellusStruct} which says that, for a special class of
noetherian local complete rings $R$ it is true that
$D(H^h_I(D(H^h_I(R))))$ is either zero or isomorphic to $R$ if
$H^l_I(R)=0$ for every $l>h=\height (I)$; the generalization is due
to Khashyarmanesh (\cite [Corollary 2.6]{khashyarmaneshPre}) and
says that $D(H^h_I(D(H^h_I(R))))\cong R$ for every noetherian local
complete ring and every ideal $I\subseteq R$ such that $H^l_I(R)=0$
for $l\neq h=\height (I)$.

We also show in section 2 that $\mu _{D(H^h_I(R))}$ is an
isomorphism if $H^l_I(R)=0$ for every $l\neq h$.

Recently there was some work on Matlis duals of local cohomology
modules (e.~g.
\cite{hellus05,hellusmArithm,hellusHabil,hellusStruct,hellusMatlisTop}).
In \cite[Corollary 1.1.4]{hellusHabil} the following was proved: If,
for some $h\in \N$, $H^l_I(R)=0$ for all $l>h$ and $\underline
x=x_1,\dots ,x_h\in I$ is an $R$-regular sequence then one has the
following equivalence
\[
\sqrt I=\sqrt{\underline xR}\iff \underline x\hbox{ is a
}D(H^h_I(R))\hbox{ regular sequence}
\]
In section 3 we extend this equivalence:

\begin{satzo}
Let $(R,m)$ be a noetherian local complete ring and $I$ an ideal of
$R$ such that $H^l_I(R)=0$ for every $l>h:=\height (I)\geq 1$; let
$\underline x=x_1,\dots ,x_h\in I$ be an $R$-regular sequence. Set
$D:=D(H^h_I(R))$. The following statements are equivalent:

(i) $\sqrt I=\sqrt {\underline xR}$; in particular, $I$ is a
set-theoretic complete intersection.

(ii) $\underline x$ is a $D$-regular sequence.

(iii) The canonical map $D/\underline xD\to H^h_{\underline xR}(D)$
(coming from $H^h_{\underline xR}(D)=\dirlim _{l\in \N}D/\underline
x^lD$) is injective.

(iv) The canonical map $\{ r\in R\vert \forall _{l\in \N}\exists
_{s\in \N}\ r\cdot I^s\subseteq \underline x^lR\} \to \Gamma
_I(R/\underline xR)$ is surjective.
\end{satzo}

The equivalence of (ii) and (iii) is inspired by a result of Marley
and Rogers (\cite[Prop. 2.3]{marleyGoren} which is a version of (ii)
$\iff $ (iii) for (arbitrary) finitely generated modules $D$.

\section{Endomorphism Rings}

\begin{definition}
(i) Let $R$ be a ring and $M$ an $R$-module. The map
\[ R\to \End _R(M), r\mapsto \hbox{ multiplication by }r\hbox{ on
}M\] is an $R$-algebra homomorphism and will be denoted by $\mu _M$.

(ii) Let $R$ be a local ring and $M$ an $R$-module. $M$ has a
canonical embedding
\[ M\to D(D(M))=D^2(M), m\mapsto (\varphi \mapsto \varphi (m))\]
into its bidual; this map will be denoted by $\iota _M$. We will
consider $M$ as a submodule of $D^2(M)$ via $\iota _M$.

(iii) Let $(R,m)$ be a noetherian local ring and $\underline
x=x_1,\dots ,x_h$ a sequence of elements of $R$. For every
$R$-module $M$ there is a canonical map
\[ M/\underline xM\buildrel \iota_{M,\underline x}\over \to H^h_{\underline xR}(M)\]
(coming from the description $H^h_{\underline xR}(M)=\dirlim _{l\in
\N }M/(x_1^l,\dots ,x_h^l)M$, where the transition maps are induced
by multiplication by $x_1\cdot \dots \cdot x_h$).
\end{definition}

\begin{satz}
\label{DHDHR}Let $(R,m)$ be a noetherian local complete ring and $I$
an ideal of $R$ such that $H^l_I(R)=0$ for every $l\neq h$ ($h$ is
necessarily the height of $I$ then). Set $H:=H^h_I(R)$.

(i) $\Hom _R(H,\iota _H):\End_ R(H)\to \Hom _R(H,D^2(H))$ is an
isomorphism.

(ii) There is a canonical isomorphism
\[ \gamma _H:\Hom _R(H,D^2(H))\to D(H^h_I(D(H)))\ \ .\]

(iii) $\mu _H:R\to \End _R(H)$ is an isomorphism of $R$-algebras.

Consequently there is a canonical isomorphism
\[ \gamma _H\circ \Hom _R(H,\iota _H)\circ \mu _H:R\to
D(H^h_I(D(H)))\ \ \]
\end{satz}
{\it Proof.} (i) It is clear that $\Hom _R(H,\iota _H)$ is
injective. To show surjectivity, let $\varphi \in \Hom _R(H,D^2(H))$
be arbitrary; let $x\in H$ be arbitrary and $n\in \N$ such that
$I^n\cdot x=0$. This implies $I^n\cdot \varphi (x)=0$, i.~e.
\[ \varphi (x)\in (0:_{D^2(H)}I^n)=D^2((0:_HI^n))=(0:_HI^n)\subseteq
H\] (the first equality follows from exactness of $D$, for the
second equality we remark that $(0:_HI^n)$ is finitely generated, as
the spectral sequence
\[ E^{p,q}_2=\Ext ^p_R(R/I^n,H^q_I(R))\Rightarrow \Ext
^{p+q}_R(R/I^n,R)\] shows $(0:_HI^n)=\Ext ^h_R(R/I^n,R)$). This
means that the image of $\varphi $ is contained in $H\subseteq
D^2(H)$, which was precisely what we had to show.

(ii) Hom-Tensor adjointness shows
\[ \Hom _R(H,D^2(H))=D(H\otimes _RD(H))\ \ .\]
On the other hand, our hypotheses imply $H^l_I=0$ for every $l>h$;
in particular $H^h_I$ is right exact, we get
\[ H\otimes _RD(H)=H^h_I(R)\otimes _RD(H)=H^h_I(D(H))\]
and statement (ii) is clear now.

(iii) \cite[Corollary 2.6]{khashyarmaneshPre} implies that there
exists an isomorphism of $R$-modules
\[ D(H^h_I(D(H)))\cong R\ \ .\]
Therefore, (i) and (ii) show that the $R$-module $\End _R(H)$ is
free of rank one. Fix any isomorphism $R\cong \End _R(H)$ and let
$\psi \in \End _R(H)$ be the element corresponding to $1\in R$. In
particular there exists a (unique) $x\in R$ such that $\id _H=x\circ
\psi $, where $x$ is multiplication by $x$ on $H$. This implies
$\varphi =x\circ \psi \circ \varphi $ for every $\varphi \in \End
_R(H)$; in particular, multiplication by $x$ is surjective on $\End
_R(H)\cong R$, $x$ is a unit in $R$ and $\mu _H$ is bijective.\hfill
$\square $

\begin{coro}
In the situation of theorem \ref{DHDHR} the endomorphism ring of
$H^h_I(R)$ is canonically isomorphic to $R$, in particular it is
commutative and $\Ann _R(H^h_I(R))=0$ holds.
\end{coro}

Let $(R,m)$ be a noetherian local ring and $M$ an $R$-module.
Consider the sequence of $R$-modules
\[ R\to \End _R(D(M))\to \End _R(D^2(M))\to \Hom _R(M,D^2(M))\ \ ,\]
where the first map is $\mu _{D(M)}$, the second is given by $\alpha
\mapsto D(\alpha )$ and the third is restriction to $M\subseteq
D^2(M)$. The composition of the second and third is always
injective:

\begin{lemma}
\label{inj}
Let $(R,m)$ be a noetherian local ring and $M$ an
$R$-module. The $R$-linear map
\[ \End _R(D(M))\to \Hom _R(M,D^2(M)), \varphi \mapsto (m\mapsto (D(\varphi
))(\iota _M (m)))\] is injective.
\end{lemma}
{\it Proof.} This is straightforward: Let $\varphi $ be in the
kernel of the above map; this means that for all $m\in M$ and for
all $\psi \in D(M)$ one has $\varphi (\psi)(m)=0$, i.~e. $\varphi
=0$.\hfill $\square $

\bigskip
We apply this injectivity in the case where the local ring $(R,m)$
is complete, $M:=H:=H^h_I(R)$ and $I$ is an ideal of $R$ such that
$H^l_I(R)=0$ for every $l\neq h$; we get $R$-linear maps
\begin{eqnarray*}
R\to \End _R(D(H))\to \End _R(D^2(H))&\to &\Hom
_R(H,D^2(H))\\
&\buildrel
\ref{DHDHR} (i)\over =&\End _R(H)\\
&\buildrel \ref{DHDHR} (iii)\over =&R
\end{eqnarray*}
The composition of all these maps is clearly $\id _R$. Thus, the
injectivity statement from lemma \ref{inj} shows:

\begin{satz}
Let $(R,m)$ be a noetherian local complete ring and $I$ an ideal of
$R$ such that $H^l_I(R)=0$ for every $l\neq h$. Then the canonical
map
\[ \mu _{D(H^h_I(R))}:R\to \End _R(D(H^h_I(R)))\]
is an isomorphism of $R$-algebras.
\end{satz}

\section{Complete Intersections and Local Cohomology}

We need a couple of lemmata and remarks before we can prove theorem
\ref{crit}, which is the main result of this section:
\begin{remark}
\label{krull} Let $(R,m)$ be a noetherian local ring, $I$ an ideal
of $R$ and $M$ an $R$-module such that \[ \Supp _R(M)\subseteq
V(I)\] (where $V(I)=\{ p\in \Spec (R)\vert p\supseteq I\} $). Let
$\hat \ $ denote $I$-adic completion. Then the natural map
\[ D(M)\to \widehat {D(M)}\]
is an isomorphism; in particular, $\bigcap _{l\in \N}I^l\cdot
D(M)=0$.
\end{remark}
{\it Proof.} We have to show that the canonical map \[ D(M)\to
\invlim _{\l\in \N}(D(M)/I^lD(M))\] is bijective; but one has
\begin{eqnarray*}
D(M)&=&D(\Gamma _I(M))\\
&=&D(\dirlim _{l\in \N} \Hom _R(R/I^l,M))\\
&=&\invlim _{\l\in \N}D(\Hom _R(R/I^l,M))\\
&=&\invlim _{l\in \N}D(M)/I^lD(M)
\end{eqnarray*}
and it is easy to see that this is the canonical map $D(M)\to
\widehat {D(M)}$.\hfill $\square $

\bigskip
Let $(R,m)$ be a noetherian local ring and $\underline x=x_1,\dots
,x_h$ a sequence of elements of $R$. Marley and Rogers have shown
(\cite[Proposition 2.3]{marleyGoren}) that, for finitely generated
$M$, $\iota_{M,\underline x}$ is injective iff $\underline x$ is an
$M$-regular sequence; in this context, note that the proof of the
following lemma is strongly based on the their proof; our additional
ingredient is remark \ref{krull}.

\begin{lemma}
\label{genMarley} Let $(R,m)$ be a noetherian local ring, $I$ an
ideal of $R$, $n,h\in \N$, $\underline x=x_1,\dots ,x_h\in I$ an
arbitrary sequence and $N$ an $R$-module. Set $H:=H^n_I(N)$ and
$D:=D(H)$. The following two statements are equivalent:

(i) For every $i=1,\dots ,h$, multiplication by $x_i$ on
$D/(x_1,\dots ,x_{i-1})D$ is injective (i.~e., $\underline x$ is a
$D$-quasiregular sequence).

(ii) $D/\underline xD\buildrel \iota _{D,\underline x}\over \to
H^h_{\underline xR}(D)$ is injective.
\end{lemma}
{\it Proof.} (i) $\Rightarrow $ (ii): The finite case is well-known,
\cite[Prop. 5.2.1]{strooker90} is a reference for the general case
(note that (ii) holds trivially if $D/\underline xD=0$).

(ii) $\Rightarrow $ (i): By induction on $h$: $h=1$: Set $x=x_1$ and
let $\alpha \in D$ be such that $x\alpha =0$. We have to show
$\alpha =0$. $\alpha $ represents an element of $\ker (\iota
_{D,\underline x})$; therefore, by assumption, $\alpha \in xD$.
Choose $\alpha _1\in D$ such that $\alpha =x\alpha _1$. We conclude
$x^2\alpha _1=0$. Again, $\alpha _1$ represents an element of $\ker
(\iota _{D,\underline x})$ and so there exists $\alpha _2\in D$ such
that $\alpha _1=x\alpha _2$. Continuing this way, we get \[ \alpha
\in \bigcap _{k\in \N}x^kD\] and then $\alpha =0$, by remark
\ref{krull}. $h>1$: First of all we prove injectivity of
\[ D/(x_1,\dots ,x_{h-1})D\buildrel \iota _{D,x_1,\dots
,x_{h-1}}\over \to H^{h-1}_{(x_1,\dots ,x_{h-1})R}(D)\ \ ;\] to do
so, let $\alpha \in \ker (\iota _{D,x_1,\dots ,x_{h-1}})$ be
arbitrary. We show $\alpha \in (x_1,\dots ,x_{h-1})D+x_h^kD$ for
every $k\in \N$ by induction on $k$: $k=0$ is trivial, we assume
$k>0$ and write $\alpha =\omega +x_h^k\beta $ for some $\omega \in
(x_1,\dots ,x_{h-1})D,\beta \in D$. By our choice of $\alpha $ there
exists $t\in \N$ such that
\[ (x_1\cdot \dots \cdot x_{h-1})^tx_h^k\beta \in (x_1^{t+1},\dots
,x_{h-1}^{t+1})D\] and hence
\[ (x_1\cdot \dots \cdot x_h)^{t+k}\beta \in (x_1^{t+k+1},\dots
,x_{h-1}^{t+k+1})D\ \ .\] But $\iota _{D,\underline x}$ is
injective, we conclude $\beta \in (x_1,\dots ,x_h)D$ and our
induction on $k$ is finished:
\[ \alpha \in \bigcap _{k\in \N}((x_1,\dots ,x_{h-1})D+x_h^kD)\ \ .\]
The $R$-module
\[ D/(x_1,\dots ,x_{h-1})D=D(\Hom _R(R/(x_1,\dots ,x_{h-1})R,H))\]
is $x_hR$-adically separated by remark \ref{krull}. This means
\[ \bigcap _{k\in \N}((x_1,\dots ,x_{h-1})D+x_h^kD)=(x_1,\dots
,x_{h-1})D\] and the stated injectivity of $\iota _{D,x_1,\dots
,x_{h-1}}$ follows. The induction hypothesis shows that $x_1,\dots
,x_{h-1}$ is $D$-quasiregular; we have to show that multiplication
by $x_h$ on $D/(x_1,\dots ,x_{h-1})D$ is injective. Let $\alpha \in
D$ be such that $x_h\alpha \in (x_1,\dots ,x_{h-1})D$. We state
\[ \forall _{k\in \N}\alpha \in (x_1,\dots ,x_{h-1})D+x_h^kD\]
and prove this statement by induction on $k$. We may assume $k>0$
and write $\alpha =\omega +x_h^k\beta $ for some $\omega \in
(x_1,\dots ,x_{h-1})D,\beta \in D$. From $x_h\alpha =x_h\omega
+x_h^{k+1}\beta $ we conclude $x_h^{k+1}\beta \in (x_1,\dots
,x_{h-1})D$. Therefore,
\[ (x_1^{k+1}\cdot \dots \cdot x_h^{k+1})\beta \in (x_1^{k+2},\dots
,x_{h-1}^{k+2})D\ \ .\] But $\iota _{D,\underline x}$ is injective
and so $\beta \in (x_1,\dots ,x_h)D$, induction on $k$ is finished:
\[ \alpha \in \bigcap _{k\in \N}((x_1,\dots
,x_{h-1})D+x_hD)\buildrel \ref{krull}\over =(x_1,\dots ,x_{h-1})D\]
(note that the last equality has been explained above in a similar
situation).

\hfill $\square $

\bigskip
Let $(R,m)$ be a noetherian local complete ring, $I$ an ideal of
$R$, $h\in \N$; assume that $\underline x=x_1,\dots ,x_h\in I$ is an
$R$-regular sequence. It follows from the Grothendieck spectral
sequence belonging to the composed functors $\Gamma _I\circ \Gamma
_{\underline xR}$ that
\[ H^h_I(R)=\Gamma _I(H^h_{\underline
xR}(R))\subseteq H^h_{\underline xR}(R)\ \ .\] By applying the
functors $D$, $H^h_{\underline xR}$ and then $D$ again, we get a
monomorphism (because $D$ is exact and $H^h_{\underline xR}$ is
right exact)
\[ D(H^h_{\underline xR}(D(H^h_I(R))))\hookrightarrow D(H^h_{\underline xR}(D(H^h_{\underline
xR}(R))))\ \ .\] Because of \cite[Corollary 2.6]{khashyarmaneshPre},
there is an isomorphism
\[ D(H^h_{\underline xR}(D(H^h_{\underline
xR}(R))))\cong R\ \ .\] Clearly, this isomorphism is unique up to a
unit of $R$ and so we may consider $D(H^h_{\underline
xR}(D(H^h_I(R))))$ as an ideal of $R$ (alternatively we use theorem
\ref{DHDHR} and have a canonical isomorphism $D(H^h_{\underline
xR}(D(H^h_{\underline xR}(R))))=R$; the resulting ideal
$J_{\underline x,I}$ is the same in both cases).
\begin{definition}
\label{idealJ} In the above situation, set
\[ J_{\underline x,I}:=D(H^h_{\underline
xR}(D(H^h_I(R))))\] and consider $J_{\underline x,I}$ as an ideal of
$R$.
\end{definition}

\begin{remark}
\label{altern} Though the definition of $J_{\underline x,I}$ is
quite abstract it also has the following concrete description:
Because of right exactness of $H^h_{\underline xR}$,
\[ D(H^h_{\underline xR}(D(H^h_I(R))))=D(H^h_{\underline xR}(R)\otimes
_RD(H^h_I(R)))\] and by Hom-Tensor adjointness, the latter module is
\[ \Hom _R(H^h_{\underline xR}(R),D^2(H^h_I(R)))\ \ .\]
Now the arguments from the proof of theorem \ref{DHDHR} (i) show
\[ \Hom _R(H^h_{\underline xR}(R),D^2(H^h_I(R)))=\Hom
_R(H^h_{\underline xR}(R),H^h_I(R))\ \ .\] But $H^h_I(R)=\Gamma
_I(H^h_{\underline xR}(R))$ and we get
\[J_{\underline x,I}=\{ \varphi \in \End _R(H^h_{\underline
xR}(R))\vert \im (\varphi )\subseteq \Gamma _I(H^h_{\underline
xR}(R))\} \ \ .\] We have $\End _R(H^h_{\underline xR}(R))\buildrel
\hbox{can}\over =R$ and thus
\[ J_{\underline x,I}=\{ r\in R\vert r\cdot H^h_{\underline
xR}(R)\subseteq \Gamma _I(H^h_{\underline xR}(R))\ \ .\] Using the
description $H^h_{\underline xR}(R)=\bigcup _{l\in \N}R/\underline
x^lR$ (where $\underline x^l=x_1^l,\dots ,x_h^l$) we conclude
\[ J_{\underline x,I}=\{ r\in R\vert \forall _{l\in \N}\exists _{s\in \N}\ r\cdot I^s\subseteq
\underline x^lR\} =\bigcap _{l\in \N}(\underline x^lR:<I>)\ \ .\]
Therefore, if we restrict the canonical map $R\to R/\underline xR$
to $J_{\underline x,I}$, we get a canonical map from $J_{\underline
x,I}$ to $\Gamma _I(R/\underline xR)$:

\begin{definition}
In the above situation, the canonical map
\[ J_{\underline x,I}\to \Gamma _I(R/\underline xR)\]
is denoted by $j_{\underline x,I}$.
\end{definition}

\end{remark}

\begin{remark}
\label{eitheror} Let $(R,m)$ be a noetherian local complete ring.
Let $I$ be an ideal of $R$, $h\in \N$ and
\[ \underline x=x_1,\dots ,x_h\subseteq I\]
an $R$-regular sequence. Then
\[ I\subseteq \sqrt {\underline xR+\Ann _R(J_{\underline x,I})}\ \
.\] In particular, if $R$ is a domain and $\sqrt{\underline
xR}\subsetneq \sqrt I$, $J_{\underline x,R}=0$
\end{remark}
{\it Proof.} We use the description of $J_{\underline x,R}$ from
remark \ref{altern}. For the first statement we have to show \[
V(\Ann _R(J_{\underline x,I}))\cap V(\underline xR)\subseteq V(I)\ \
,\] i.~e. for every $r\in J_{\underline x,R}$ we have to show \[
V(\Ann _R(r))\cap V(\underline xR)\subseteq V(I)\ \ :\] Let $r\in
J_{\underline x,I}$ be arbitrary; by remark \ref{altern},
$J_{\underline x,I}=\{ r\in R\vert \forall _{l\in \N}\exists _{s\in
\N}\ r\cdot I^s\subseteq \underline x^lR\} $. For every $p\in
V(\underline xR)\setminus V(I)$ we get
\[ r\cdot R_p\subseteq \bigcap _{l\in \N}\underline x^lR_p\subseteq \bigcap _{l\in
\N}p^lR_p=0\ \ ,\] i.~e. $\Ann _R(r)\not\subseteq p$ and the first
statement is proven. The second statement follows immediately from
the first. \hfill $\square $

\begin{satz}
\label{crit} Let $(R,m)$ be a noetherian local complete ring and $I$
an ideal of $R$ such that $H^l_I(R)=0$ for every $l>h:=\height
(I)\geq 1$; let $\underline x=x_1,\dots ,x_h\in I$ be an $R$-regular
sequence (clearly, this implies $H^l_I(R)=0$ for every $l\neq h$) .
Set $D:=D(H^h_I(R))$. The following statements are equivalent:

(i) $\sqrt I=\sqrt {\underline xR}$; in particular, $I$ is a
set-theoretic complete intersection.

(ii) $\underline x$ is a $D$-regular sequence.

(iii) $D/\underline xD\buildrel \iota _{D,\underline x}\over \to
H^h_{\underline xR}(D)$ is injective.

(iv) $j_{\underline x,I}$ is surjective.

(v) $J_{\underline x,I}=R$.
\end{satz}
{\it Proof.} (i) $\iff $ (ii) was shown (for more general $R$) in
\cite[Cor. 1.1.4]{hellusHabil}.

(ii) $\iff  $ (iii): Is a special case of Lemma \ref{genMarley}
(note that $D/\underline xD=D(\Hom _R(R/\underline xR,H^h_I(R)))\neq
0$).

(iii) $\iff $ (iv): By definition, $D(H^h_{\underline
xR}(D))=J_{\underline x,I}$. We have
\[ D(D/\underline xD)=D(D(\Hom
_R(R/\underline xR,H^h_I(R))))\ \ .\] But $H^h_I(R)=\Gamma
_I(H^h_{\underline xR}(R))$ and, therefore,
\begin{eqnarray*}
\Hom _R(R/\underline xR,H^h_I(R))&=&\Gamma _I(\Hom _R(R/\underline
x,H^h_{\underline xR}(R)))\\
&=&\Gamma _I(\Ext ^h_R(R/\underline xR,R))\\
&=&\Gamma _I(R/\underline xR)
\end{eqnarray*}
(for the second and the third equality use the fact that $\underline
x$ is an $R$-regular sequence). The latter module is finitely
generated, we get
\[ D(D(\Hom _R(R/\underline xR,H^h_I(R))))=\Gamma _I(R/\underline xR)\ \
.\] Thus $D(\iota _{D,\underline x})$ is a map $J_{\underline
x,I}\to \Gamma _I(R/\underline xR)$; it is straightforward to see
that it is in fact $j_{\underline x,I}$ (to do so one should start
with the description $H^h_{\underline xR}(D)=\dirlim _{l\in
\N}(D/\underline x^lD$).\hfill $\square $
\begin{question}
In the situation of definition \ref{idealJ}, when exactly is
$J_{\underline x,I}=0$?
\end{question}

{\sc Universit\"at Leipzig, Fakult\"at f\"ur Mathematik und
Informatik, Ma\-the\-ma\-tisches Institut, Augustusplatz 10/11,
D-04109 Leipzig}

{\it E-mail}:

hellus@math.uni-leipzig.de

stueckrad@math.uni-leipzig.de
\end{document}